\documentclass[journal]{IEEEtran}

\usepackage{amsmath,amsfonts,amssymb,amsmath,latexsym,multirow}

\newtheorem{theorem}{Theorem}
\newtheorem{lemma}[theorem]{Lemma}
\newtheorem{corollary}[theorem]{Corollary}
\newtheorem{proposition}[theorem]{Proposition}
\begin{document}
\title{ Codes Associated with Special Linear Groups \qquad
         and Power Moments of Multi-dimensional Kloosterman Sums}
\author{Dae San Kim,~\IEEEmembership{Member,~IEEE}
\thanks{This work was supported by grant No. R01-2008-000-11176-0 from the Basic Research Program of the Korea Science and Engineering Foundation.}% <-this % stops a space
\thanks{The author is with the Department of Mathematics, Sogang University, Seoul 121-742, Korea(e-mail: dskim@sogang.ac.kr). }
        }
\maketitle

\begin{abstract}
In this paper, we construct the binary linear codes $C(SL(n,q))$
 associated with finite special linear groups $SL(n,q)$, with both
\emph{n,q} powers of two. Then, via Pless power moment identity and
utilizing our previous result on the explicit expression of the
Gauss sum for $SL(n,q)$, we obtain a recursive formula for the power
moments of multi-dimensional Kloosterman sums in terms of the
frequencies of weights in $C(SL(n,q))$. In particular, when $n=2$,
this gives a recursive formula for the power moments of Kloosterman
sums. We illustrate our results with  some examples.
\end{abstract}



% Will be ignored for other modes.
\IEEEpeerreviewmaketitle

%%%%%%%%%%%%%%%%%%%%%%%%%%%%%%%%%%%%%%%%%%%%%%%%%%%%%%%%%%%%%%%%%%%%
\section{Introduction}
%%%%%%%%%%%%%%%%%%%%%%%%%%%%%%%%%%%%%%%%%%%%%%%%%%%%%%%%%%%%%%%%%%%%
Let $\psi$ be a nontrivial additive character of the finite field
$\mathbb{F}_{q}$ with $q=p^{r}$ elements ( $p$ a prime), and let $m$
be a positive integer. Then the $m$-dimensional  Kloosterman sum
$K_{m}(\psi;a)$(\cite{RH}) is defined by

\begin{align*}
K_{m}(\psi;a)=\sum_{\alpha_{1},\cdots,\alpha_{m} \in
\mathbb{F}_{q}^{*}}\psi(\alpha_{1}+\cdots+\alpha_{m}+a\alpha_{1}^{-1}\cdots\alpha_{m}^{-1})\\
(a \in \mathbb{F}_{q}^{*}).
\end{align*}

In particular, if $m=1$, then $K_{1}(\psi;a)$ is simply denoted by
$K(\psi;a)$, and is called the Kloosterman sum. The Kloosterman sum
was introduced in 1926 \cite{HDK} to give an estimate for the
Fourier coefficients of modular forms.

For each nonnegative integer $h$, by $MK_{m}(\psi)^{h}$ we will
denote the $h$-th moment of the $m$-dimensional Kloosterman sum
$K_{m}(\psi;a)$. Namely, it is given by
\begin{equation*}
 MK_{m}(\psi)^{h}=\sum_{a \in \mathbb{F}_{q}^{*}}K_{m}(\psi;a)^{h}
 \end{equation*}

If $\psi=\lambda$ is the canonical additive character of
$\mathbb{F}_{q}$, then $MK_{m}(\lambda)^{h}$ will be simply denoted
by $MK_{m}^{h}$. If further $m=1$, for brevity $MK_{1}^{h}$ will be
indicated by $MK^{h}$. The power moments of Kloosterman sums can be
used, for example, to give an estimate for the Kloosterman sums and
have also been studied to solve a variety of problems in coding
theory over finite fields of characteristic two.

If $q=p$ is an odd prime, for $h\leq 4$, $MK^{h}$ was evaluated by
Sali\'{e} \cite{HS}. For details about these, the reader is referred
to Section IV.

From now on, let us assume that $q=2^{r}$. Carlitz \cite{L1}
evaluated $MK^{h}$, for $h\leq 4$, while Moisio computed $MK^{6}$ in
 \cite{MK}. Recently, Moisio was able to find explicit expressions of
$MK^{h}$, for $h \leq 10$ (cf. \cite{M1}). This was done, via Pless
power moment identity, by connecting moments of Kloosterman sums and
the frequencies of weights in the binary Zetterberg code of length
$q+1$, which were known by the work of Schoof and Vlugt in
\cite{RM}.

In this paper, we adopt Moisio's idea to show the following theorem
giving a recursive formula for the power moments of
multi-dimensional Kloosterman sums. To do that, we construct the
binary linear code $C(SL(n,q))$ associated with the special linear
group $SL(n,q)$, and express those power moments in terms of the
frequencies of weights in the code. Here, in addition to the
assumption $q=2^{r}$, we restrict $n$ to be $n=2^{s}$. Then, thanks
to our previous result on the explicit expression of ``Gauss sum''
for the special linear group \cite{DS}, we can express the weight of
each codeword in the dual $C^{\perp}(SL(n,q))$ of $C(SL(n,q))$, in
terms of $(n-1)$-dimensional Kloosterman sums. Then our formula
follows immediately from the Pless power
moment identity.\\

\begin{theorem}\label{A}
Let $n=2^{s}$, $q=2^{r}$. Then, for all positive integers $h$, we
have the following recursive formula for the moments of
multi-dimensional Kloosterman sums $MK_{n-1}^{h}$:
\begin{align}\label{a}
  \begin{split}
q^{\binom{n}{2}h}&MK_{n-1}^{h}=\sum_{i=0}^{h-1}(-1)^{h+i+1}{\binom{h}{i}}N^{h-i}q^{\binom{n}{2}i}MK_{n-1} ^{i}\\
                 &+q\sum_{i=0}^{min\{N,h\}}(-1)^{h+i}C_{i}\sum_{t=i}^{h}t!S(h,t)2^{h-t}{\binom{N-i}{N-t}}.
  \end{split}
\end{align}

Here $N=q^{\binom{n}{2}}\prod_{j=2}^{n}(q^{j}-1)$ is the order of
$SL(n,q)$, and $S(h,t)$ indicates  the Stirling number of the second
kind given by

\begin{equation}\label{b}
S(h,t)=\frac{1}{t!}\sum_{j=0}^{t}(-1)^{t-j}{\binom{t}{j}}j ^{h} .
\end{equation}

In addition, $\{C_{i}\}_{i=0}^{N}$ denotes  the weight distribution
of the code $C=C(SL(n,q))$, which is given by

\begin{equation*}
C_{i}=\sum_{}^{}\prod_{\beta \in \mathbb{F}_{q}}{\binom
{n_{\beta}}{\nu
_{\beta }}}(0 \leq i \leq N),\\
\end{equation*}

where the sum runs over all the sets of nonnegative integers
$\{\nu_{\beta}\}_{\beta \in \mathbb{F}_{q}}$ satisfying $\sum_{\beta
\in \mathbb{F}_{q}} ^{}\nu_{\beta }=i$ \, and \, $\sum_{\beta \in
\mathbb{F}_{q}}^{}\nu_{\beta}\beta =0$ \; (an identity in
$\mathbb{F}_{q}$), and

\begin{align*}
\begin{split}
n_{\beta}&=|\{g \in SL(n,q)|tr(g)=\beta\}|\\
&=q^{{\binom{n}{2}}-1}\{\prod_{j=2}^{n}(q^{j} -1)+1+q\theta(
\beta)\},\\
\end{split}
\end{align*}
with\begin{equation*} \theta(\beta)=
\begin{cases}
 K_{n-2}(\lambda;\beta ^{-1}),& \beta \neq 0,\\
0,& \beta =0.
\end{cases}
\end{equation*}

Here we understand that
$K_{0}(\lambda;\beta^{-1})=\lambda(\beta^{-1})$. In addition, from
now on we agree that ${\binom{b}{a}}=0$, if $b<a$ .\\
\end{theorem}

\begin{corollary}\label{B}
Let $q=2^{r}$. Then, for all positive integers $h$ , we have the
following recursive formula for the moments of  Kloosterman sums
$MK^{h}$:

\begin{align}\label{c}
\begin{split}
 q^{h}&MK^{h}=\sum_{i=0}^{h-1}(-1)^{h+i+1}{\binom{h}{i}}N^{h-i}q^{i}MK^{i}\\
 &+q\sum_{i=0}^{min\{N,\,h\}}(-1)^{h+i}C_{i}\sum_{t=i}^{h}t!S(h,t)2^{h-t}{\binom
{N-i}{N-t}}.
\end{split}
\end{align}

Here $N=q(q^2-1)$ is the order of $SL(2,q)$, $S(h,t)$ indicates the
Stirling number of the second kind as in ({\ref{b}}), and
$\{C_{i}\}_{i=0}^{N}$ denotes the weight distribution of the code
$C=C(SL(2,q))$, which is given by

\begin{align*}
\begin{split}
C_{i}=\sum_{}^{}{\binom{q^{2}}{\nu_{0}}}\prod_{tr(\beta^{-1})=0}{\binom{q^{2}+q}{\nu_{\beta}}}
\prod_{tr(\beta^{-1})=1}{\binom{q^{2}-q}{\nu_{\beta }}}\\
 (0 \leq i \leq N),
\end{split}
\end{align*}

where the sum runs over all the sets of nonnegative integers
$\{\nu_{\beta}\}_{\beta \in \mathbb{F}_{q}}$ satisfying $\sum_{\beta
\in \mathbb{F}_{q}} ^{}\nu_{\beta }=i$ and $\sum_{\beta \in
\mathbb{F}_{q}}^{}\nu_{\beta}\beta =0$, and the first and second
product run respectively over the elements $\beta \in
\mathbb{F}_{q}^{*}$, with $tr(\beta^{-1})=0$ and $tr(\beta^{-1})=1$.
\end{corollary}

%%%%%%%%%%%%%%%%%%%%%%%%%%%%%%%%%%%%%%%%%%%%%%%%%%%%%%%%%%%%%%%%%%%%
\section{Preliminaries}
%%%%%%%%%%%%%%%%%%%%%%%%%%%%%%%%%%%%%%%%%%%%%%%%%%%%%%%%%%%%%%%%%%%%
The following notations will be used throughout this paper except in
Section IV, where $q$ is allowed to be any prime powers.\\

$n=2^s (s \in \mathbb{Z}_{>0})$,\\

$q=2^r (r \in \mathbb{Z}_{>0})$,\\

$SL(n,q)$= the special linear group,\\

$N$=$q^{{\binom{n}{2}}} \prod_{j=2}^{n}(q ^{j}-1)$ the order of
$SL(n,q)$,\\

$Tr(g)$= the matrix trace for $g \in SL(n,q)$,\\

$tr(x)=x+x^2+\cdots+x^{2^{r-1}}$ the trace function
$\mathbb{F}_{q}\longrightarrow \mathbb{F}_{2}$,\\

$\lambda(x)$=$(-1)^{tr(x)}$ the canonical additive character of
$\mathbb{F}_{q}$.\\

Let $g_{1}$, $g_{2}$,$\cdots$, $g_{N}$ be a fixed ordering of the
elements in $SL(n,q)$. Let $C=C(SL(n,q))$ be the binary linear code
of length $N$, defined by:\\
\begin{equation}\label{d}
C=C(SL(n,q))= \{u \in \mathbb{F}_{2}^{N}|u \cdot v=0\},
\end{equation}
where
\begin{equation} \label{e} v=(Tr(g_{1}),Tr(g_{2}),\cdots
,Tr(g_{N}))\in \mathbb{F}_{q} ^{N}.
\end{equation}
\smallskip

\begin{theorem}[Delsarte, \cite{FN}]\label{C} Let $B$ be a linear code over
 $\mathbb{F}_{q}$. Then $(B|_{\mathbb{F}_{2}})^{\bot}=tr(B^{\bot})$.\\
\end{theorem}

From Delsarte's theorem, the next result follows immediately.\\

\begin{proposition}\label{D}
The dual $C^{\bot}=C^{\bot}(SL(n,q))$ of $C=C(SL(n,q))$ is given by
\begin{align*}
\begin{split}
 C^{\bot}=\{c(a)=(tr(aTr(g_{1}))&, tr(aTr(g _{2})), \cdots , \\
                                 &tr(aTr(g_{N})))| \,\; a \in \mathbb{F}_{q}\}.\\
\end{split}
\end{align*}

The next Proposition is stated in Theorem 6.1 of {\cite{DS}}. But we
slightly modified the expression there.\\
\end{proposition}

\begin{proposition}\label{E}
Let $n_{\beta}=|\{g \in SL(n,q)| \; Tr(g)= \beta \}|$, \qquad for
each $\beta \in \mathbb{F}_{q}$. Then
\begin{equation*}
n_{\beta}=q^{\binom{n} {2}-1}\{\prod_{j=2}^{n}(q^{j}
-1)-(q-1)^{n-1}+q \delta(n-1,q;\beta )\},
\end{equation*}

where
\begin{align*}
\begin{split}
\delta(n-1,q;\beta )=|\{( \alpha _{1}, \cdots,&\alpha_{n-1})\in (\mathbb{F}_{q}^{*})^{n-1}|\alpha_{1}+\cdots\\
&+\alpha_{n-1}+\alpha_{1}^{-1} \cdots \alpha_{n-1}^{-1}=\beta \}|.
\\ \\
\end{split}
\end{align*}
\end{proposition}

The following corollary is immediate from Proposition 5.\\

\begin{corollary}\label{F}
The map $Tr$ : $SL(n,q)\longrightarrow \mathbb{F}_{q}$ given by $g
\longmapsto Tr(g)$
is surjective.\\
\end{corollary}

\begin{proposition}\label{G}
The map $\mathbb{F}_{q} \longrightarrow C^{\bot}(SL(n,q))$ given by
$a \longmapsto c(a)$ is an $\mathbb{F}_2$-linear isomorphism.
\end{proposition}
\begin{proof}
It is $\mathbb{F}_2$-linear and surjective. Let $a$ be in the kernel
of the map. Then $tr(aTrg)=0$, for all $g \in SL(n,q)$. In view of
Corollary \ref{F}, $tr(a \alpha)=0$, for all $\alpha \in
\mathbb{F}_{q}$. As $tr: \mathbb{F}_{q} \longrightarrow \mathbb{F}_{2}$ is surjective, $a=0$.\\
\end{proof}

The next theorem is about the Gauss sum for $SL(n,q)$ , and is one
of the main results of the paper {\cite{DS}}.\\

\begin{theorem}\label{H}
Let $\psi$ be any nontrivial additive character of $\mathbb{F}_{q}$.
Then
\begin{equation*}
\sum_{g \in SL(n,q)}^{} \psi (Tr(g))=q^{\binom {n}{2}}
K_{n-1}(\psi;1).
\end{equation*}
For the following lemma, observe that $(n,q-1)=1$.\\
\end{theorem}

\begin{lemma}\label{I}
The map $a\longmapsto a^n : \mathbb{F}_{q}^{*} \longrightarrow
\mathbb{F}_{q}^{*}$ is a bijection.\\
\end{lemma}

For the proof of the next proposition and the following, we borrowed
an idea from the proof of Theorem 6.1 in {\cite{M3}}.\\

\begin{proposition}\label{J}
For $a \in \mathbb{F}_{q}^{*}$, the Hamming weight of the codeword
\begin{equation*}
c(a)=(tr(aTr(g_{1})),tr(aTr(g_{2})),\cdots,tr(aTr(g_{N})))
\end{equation*}
is given by(cf. Proposition \ref{D}):
\begin{equation}\label{f}
w(c(a))=\frac{1}{2}(N-q^{\binom{n}{2}}K_{n-1}(\lambda ; a)).
\end{equation}
\end{proposition}
\begin{proof}
\begin{align*}
\begin{split}
&w(c(a))=\frac{1}{2}\sum_{i=1}^{N}(1-(-1)^{tr(aTr(g_i))}) \qquad \qquad \\
        &=\frac{N}{2}-\frac{1}{2} \sum_{g \in SL(n,q)}\lambda(aTr(g)) \qquad \qquad \\
        &=\frac{N}{2}-\frac{1}{2}q^{\binom{n}{2}}K_{n-1}(\psi ;1)\qquad \qquad \\
\end{split}
\end{align*}
\qquad \qquad \qquad (Theorem \ref{H}, with $\psi(x)$=$\lambda(ax)$)
\begin{align*}
\begin{split}
 =\frac{N}{2}-\frac{1}{2}q^{\binom{n}{2}} \sum_{x_{1}, \cdots, x_{n-1} \in \mathbb{F}_{q}^{*}}^{} &\lambda(ax_{1}+ \cdots ax_{n-1}\\
                                                                                                           &+ax_{1}^{-1} \cdots x_{n-1}^{-1})\\
 =\frac{N}{2}-\frac{1}{2}q^{\binom{n}{2}} \sum_{x_{1}, \cdots, x_{n-1} \in \mathbb{F}_{q}^{*}}^{} &\lambda(x_{1}+ \cdots+ x_{n-1}\\
                                                                                                           &+a^{n}x_{1}^{-1} \cdots x_{n-1}^{-1})\\
\end{split}
\end{align*}
\begin{align}\label{g}
\begin{split}
\qquad
  =\frac{N}{2}-\frac{1}{2}q^{\binom{n}{2}}\sum_{x_1,\cdots,x_{n-1} \in \mathbb{F}_{q}^{*}}& \lambda (x_1^n+\cdots+x_{n-1}^n\\
                                                                              &+a^{n}x_1^{-n} \cdots x_{n-1}^{-n})\\
\end{split}
\end{align}
\qquad \qquad \qquad \qquad \qquad \qquad (by Lemma \ref{I})
\begin{align}\label{h}
\begin{split}
 =\frac{N}{2}-\frac{1}{2}q^{\binom{n}{2}} \sum _{x_1,\cdots,x_{n-1} \in \mathbb{F}_{q}^{*}}& \lambda((x_1+ \cdots+x_{n-1}\\
                                                                                &+ax_1^{-1}\cdots x_{n-1}^{-1})^n)\\
\end{split}
\end{align}

\begin{align}\label{i}
\begin{split}
 =\frac{N}{2}-\frac{1}{2}q^{\binom {n}{2}} \sum_{x_1,\cdots,x_{n-1} \in \mathbb{F}_{q}^{*}}& \lambda(x_1+ \cdots+x_{n-1}\\
                                                                                &+ax_1^{-1}\cdots x_{n-1}^{-1})\\
\end{split}
\end{align}
\qquad \qquad \qquad \qquad \qquad({\cite{RH}}, Theorem 2.23(v))
\begin{equation*}
 =\frac{N}{2}-\frac{1}{2}q^{\binom {n}{2}}K_{n-1}(\lambda; a). \qquad \qquad \qquad \qquad
\end{equation*}
\end{proof}

We are ready to determine $ \delta(n-1,q; \beta)$, which appears in
Proposition \ref{E}.\\
\\
\\
\begin{proposition}\label{K}
 For each $ \beta \in \mathbb{F}_{q}$, let
\begin{align*}
\begin{split}
\delta(n-1,q; \beta )=|\{(\alpha_{1},\cdots,&\alpha_{n-1}) \in (\mathbb{F}_{q}^{*})^{n-1}|\alpha_{1}+\cdots \\
                                            &+\alpha_{n-1}+ \alpha_{1}^{-1} \cdots \alpha_{n-1}^{-1}= \beta \} |.\\
\end{split}
\end{align*}
Then \begin{equation*} \delta(n-1,q;0)=q^{-1}\{(q-1)^{n-1}+1 \},
\end{equation*}
and, for $ \beta \in \mathbb{F}_{q}^{*}$,
\begin{equation*}
\delta(n-1,q;\beta)=K_{n-2}(\lambda ;\beta^{-1})+q^{-1}
\{(q-1)^{{}^{n-1}}+1 \},
\end{equation*}
where $ K_{0}(\lambda ; \beta^{-1})$ = $\lambda( \beta^{-1})$ by
convention. \\
\end{proposition}

\begin{proof}
\begin{align*}
\begin{split}
&q \delta(n-1,q; \beta)\\
 &=\sum_{\alpha_{1},\cdots,\alpha_{n-1} \in \mathbb{F}_{q}^{*}}^{}\sum_{\alpha \in \mathbb{F}_{q}}^{} \lambda(\alpha(\alpha_{1}+\cdots+\alpha_{n-1}\\
 & \qquad \qquad \qquad \qquad \qquad \qquad+\alpha_{1}^{-1}\cdots \alpha_{n-1}^{-1}-\beta )) \\
 &=\sum_{\alpha \in \mathbb{F}_{q}} \lambda(-\alpha \beta)\sum_{\alpha_{1},\cdots, \alpha_{n-1} \in \mathbb{F}_q^{*} }\lambda(\alpha \alpha_{1}+\cdots+\alpha \alpha_{n-1}\\
 & \qquad \qquad \qquad \qquad \qquad \qquad \qquad +\alpha \alpha_1^{-1} \cdots \alpha_{n-1}^{-1}) \\
 &=\sum_{\alpha \in \mathbb{F}_q^{*}}\lambda(- \alpha \beta)\sum_{\alpha_{1}, \cdots, \alpha_{n-1} \in \mathbb{F}_q^{*} } \lambda(\alpha_{1}+\cdots+\alpha_{n-1}\\
 & \qquad \qquad \qquad \qquad \qquad +\alpha^{n}\alpha_1^{-1} \cdots \alpha_{n-1}^{-1})+(q-1)^{n-1} \\
 &=\sum_{\alpha \in \mathbb{F}_q^{*}}\lambda(- \alpha \beta)\sum_{\alpha_{1},\cdots, \alpha_{n-1} \in \mathbb{F}_{q}^{*} }\lambda(\alpha_{1}+\cdots+\alpha_{n-1}\\
 & \qquad \qquad \qquad \qquad \qquad +\alpha \alpha_1^{-1}\cdots \alpha_{n-1}^{-1})+(q-1)^{n-1} \\
\end{split}
\end{align*}
\qquad \qquad (following the steps in (\ref{g})-(\ref{i}))

\begin{align*}
\begin{split}
 =&\sum_{\alpha_{1},\cdots,\alpha_{n-1} \in \mathbb{F}_{q}^{*}}^{} \lambda(\alpha_{1}+\cdots+\alpha_{n-1})\\
   & \qquad \times \sum_{\alpha \in \mathbb{F}_{q}^{*}} ^{} \lambda(\alpha(\alpha_{1}^{-1}\cdots \alpha_{n-1}^{-1}- \beta))+(q-1)^{n-1}\\
 =&\sum_{\alpha_{1},\cdots,\alpha_{n-1} \in \mathbb{F}_{q}^{*}}^{} \lambda(\alpha_{1}+\cdots+\alpha_{n-1})\\
  &  \qquad \times \sum_{\alpha  \in \mathbb{F}_{q}}^{} \lambda(\alpha(\alpha_{1}^{-1}\cdots \alpha_{n-1}^{-1}-\beta))\\
  & \qquad - \sum_{\alpha_{1},\cdots,\alpha_{n-1} \in \mathbb{F}_{q}^{*}}^{} \lambda(\alpha_{1}+\cdots+\alpha_{n-1})+(q-1)^{n-1}\\
\end{split}
\end{align*}
\begin{equation}\label{j}
 =q \sum \lambda(\alpha_{1}+\cdots+ \alpha_{n-1})+1+(q-1)^{n-1} \qquad \qquad
\end{equation}

The sum in (\ref{j}) runs over all $ \alpha_{1},\cdots,\alpha_{n-1}
\in \mathbb{F}_{q}^{*}$ satisfying $\alpha_{1}^{-1}\cdots
\alpha_{n-1}^{-1}$ = $\beta$, so that it is given by

\begin{equation*}
\begin{cases}
 0, & \text{if $\beta=0$},\\
 K_{n-2}(\lambda ; \beta^{-1}),& \text {if $\beta \neq 0$, and $n>2$,}\\
 \lambda(\beta^{-1}),& \text{if $\beta \neq0$, and $n=2$.}
\end{cases}
\end{equation*}
So we get the desired result.\\
\end{proof}

Combining Propositions \ref{E} and \ref{K}, we get the following
corollary. \\

\begin{corollary}\label{L}
Let
\begin{equation*}
n_{\beta }=|\{ g \in SL(n,q)| \;Tr(g)= \beta \}|,
\end{equation*}

for each $\beta \in \mathbb{F}_{q}$. Then

\begin{equation}\label{k}
n_{\beta}=q^{\binom{n}{2}-1}\{\prod_{j=2}^{n}(q^{j}-1)+1+q
\theta(\beta)\},
\end{equation}
where
\begin{equation*} \theta(\beta)=
\begin{cases}
 K_{n-2}(\lambda; \beta ^{-1}),& \beta \neq 0,\\
0,& \beta =0,
\end{cases}
\end{equation*}
with the convention that $K_{0}(\lambda; \beta^{-1}) =
\lambda(\beta^{-1})$.
\end{corollary}
%%%%%%%%%%%%%%%%%%%%%%%%%%%%%%%%%%%%%%%%%%%%%%%%%%%%%%%%%%%%%%%%%%%%
\section{Proof of Main Results}
%%%%%%%%%%%%%%%%%%%%%%%%%%%%%%%%%%%%%%%%%%%%%%%%%%%%%%%%%%%%%%%%%%%%
In this section, we will derive the recursive formula ({\ref{a}})
for the power moments of multi-dimensional Kloosterman sums which is
expressed in terms of the frequencies $C_{i}$ of weights in the code
$C=C(SL(n,q))$.\\

\begin{theorem}[Pless power moment identity, {\cite{FN}}]\label{M}
Let $B$ be an $q$-ary $[n,k]$ code, and let $B_{i}$(resp.
$B_{i}^{\bot})$ denote the number of codewords of weight $i$ in $B$
(resp. in $B^{\bot}$). Then, for $h=0,1,2,\cdots$,
\begin{align}\label{l}
\begin{split}
\sum_{i=0}^{n}i^{h}B_{i}=\sum_{i=0}^{min
\{n,h\}}(-1)^{i}B_{i}^{\bot}\sum_{t=i}^{h}&t!S(h,t)q^{k-t}\\
                                          &\times (q-1)^{t-i} {\binom {n-i}{n-t}},\\
\end{split}
\end{align}
where $S(h,t)$ is the Stirling number of the second kind defined in
({\ref{b}}).\\
\end{theorem}

\begin{theorem}[\cite{GJ}]\label{N}
Let $q=2^{r}$, with $r \geq 2$. Then the range $R$ of
$K(\lambda;a)$, as $a$ varies over $\mathbb{F}_{q}^{*}$, is given by
\begin{equation*}
R=\{t \in \mathbb{Z} | \; |t |<2 \sqrt {q }, \; t \equiv -1
(\textmd{mod} \; 4)\}.
\end{equation*}
In addition, each value $t \in R$ is attained exactly $H(t^2 -q)$
times, where $H(d)$ is the Kronecker class number of $d$.\\
\end{theorem}

\begin{theorem}[\cite{L2}]\label{O} For the canonical additive character
$\lambda $ of $\mathbb{F}_{q}$, and $a \in \mathbb{F}_{q}^{*}$,
\begin{equation}\label{m}
K_{2}(\lambda ;a)=K(\lambda ;a)^{2}-q.
\end{equation}

Let $u$ = $(u_{1},\cdots,u_{N}) \in \mathbb{F}_{2}^{N}$, with
$\nu_{\beta}$ 1's in the coordinate places where $Tr(g_{j})=\beta$,
for each $\beta \in \mathbb{F}_{q}$. Then we see from the definition
of the code $C=C(SL(n,q))$(cf. ({\ref{d}}),({\ref{e}})) that $u$ is
a codeword with weight $i$ if and only if $\sum_{\beta \in
\mathbb{F}_{q}}^{} \nu_{\beta}=i$ and $\sum_{\beta \in
\mathbb{F}_{q}}^{} \nu_{\beta}\beta=0$ (an identity in
$\mathbb{F}_{q}$). As there are $\prod_{\beta \in
\mathbb{F}_{q}}{\binom {n_{\beta}}{\nu_{\beta}}}$ many such
codewords with
weight $i$, we obtain the following theorem.\\
\end{theorem}

\begin{theorem}\label{P}
Let $\{C_{i}\}_{i=0}^{N}$ be the weight distribution of the code $
C=C(SL(n,q))$. Then, for $0 \leq i \leq N$,
\begin{equation}\label{n}
C_{i}=\sum_{}^{} \prod_{\beta \in
\mathbb{F}_{q}}{\binom{n_{\beta}}{\nu_{\beta}}},
\end{equation}
where $n_{\beta}$ is as in (\ref{k}), and the sum runs over all the
sets of nonnegative integers $\{\nu_{\beta}\}_{\beta \in
\mathbb{F}_{q}}$ satisfying
\begin{equation}\label{w}
 \sum_{\beta \in \mathbb{F}_{q}}^{} \nu_{\beta}=i \; \textrm{and}
 \sum_{\beta \in \mathbb{F}_{q}}^{}\nu_{\beta} \beta=0(\textrm{an}\; \textrm{identity} \; \textrm{in}\; \mathbb{F}_{q}).
\end{equation}\\
\end{theorem}

\begin{corollary}\label{W}
Let $\{C_{i}\}_{i=0}^{N}$ be the weight distribution of the code
$C=C(SL(n,q))$. Then, for $0 \leq i \leq N$, $C_{i}=C_{N-i}$.\\
\end{corollary}

\begin{proof}
Under the replacements $\nu_{\beta} \rightarrow
n_{\beta}-\nu_{\beta}$, for all $\beta \in \mathbb{F}_{q}$, the
first sum in (\ref{w}) is changed to $N-i$, while the second one in
(\ref{w}) and the summands in (\ref{n}) are left unchanged. Here the
second sum in (\ref{w}) is left unchanged, since $\sum_{\beta \in
\mathbb{F}_{q}}^{} n_{\beta} \beta =0$, as one can see by using the
explicit expression of $n_{\beta}$ in (\ref{k}).\\
\end{proof}
\begin{corollary}\label{Q}
Let $\{C_{i}\}_{i=0}^{N}$ be the weight distribution of the code
$C=C(SL(2,q))$. Then, for $0 \leq i \leq N$,

\begin{equation*}
C_{i}=\sum_{}^{}{\binom {q^{2}}{\nu_{0}}} \prod_{tr(\beta ^{-1}
)=0}{\binom {q^{2}+q}{\nu_{\beta}}} \prod_{tr(\beta^{-1})=1}{\binom
{q^{2}-q}{\nu_{\beta}}},
\end{equation*}
where the sum runs over all the sets of nonnegative integers $\{\nu
_{\beta}\}_{\beta \in \mathbb{F}_{q}}$ satisfying $\sum_{\beta \in
\mathbb{F}_{q}}^{}\nu_{\beta}=i$ and $\sum_{\beta \in
\mathbb{F}_{q}}^{}\nu_{\beta } \beta=0$, and the first and second
product run respectively over the elements $\beta \in
\mathbb{F}_{q}^{*}$, with $tr(\beta^{-1} )=0$ and
$tr(\beta^{-1} )=1$.\\
\end{corollary}

\begin{proof}
For $n=2$, we see from (\ref{k}) that $n_{\beta}$ is given by
\begin{equation*}
n_{\beta}=
\begin{cases}
q^{2} ,& \text {if $\beta=0$,}\\
q^{2}+q ,& \text {if $tr(\beta^{-1})=0$,}\\
q^{2}-q ,& \text{if $tr(\beta^{-1})=1$.}
\end{cases}
\end{equation*}
\end{proof}

\begin{corollary}\label{R}
Assume that $ r \geq 2$, and that $\{C_{i}\}_{i=0}^{N}$ is the
weight distribution of the code $C=C(SL(4,q))$. Then, for $ 0 \leq i
\leq N$,
\begin{equation}\label{o}
C_{i}=\sum {\binom {m_{0}}{\nu_{0}}} \prod_{\substack{|t |<
2\sqrt{q}\\t \equiv -1 (4)}} \qquad \prod_{K(\lambda ;
\beta^{-1})=t} {\binom{m_{t}}{\nu_{\beta}}},
\end{equation}
where the sum runs over all the sets of nonnegative integers $\{\nu
_{\beta}\}_{\beta \in \mathbb{F}_{q}}$ satisfying $\sum_{\beta \in
\mathbb{F}_{q}}^{} \nu_{\beta}=i$ and $\sum_{\beta \in
\mathbb{F}_{q}}^{}\nu_{\beta}\beta=0$,
\begin{equation*}
m_{0}=n_{0}=q^{5}\{\prod_{j=2}^4(q^j-1)+1\},
\end{equation*}
and
\begin{equation*}
m_{t}=q^{6} \{q^{2}(q^{2}-1)(q^{4}-q-1)+t^{2}\},
\end{equation*}
for all integers $t$ satisfying $|t|<2 \sqrt {q}$ and
$t\equiv-1(4)$.\\

\begin{proof}
Note here that, for $n=4$,  and $\beta \in \mathbb{F}_{q }^{*}$,
\begin{equation*}
n_{\beta}=q^{5}\{\prod_{j=2}^{4}(q^{j}-1)+1+qK_{2}(\lambda;
\beta^{-1})\} \qquad \qquad \qquad
\end{equation*}
\begin{equation}\label{p}
=q^{5}\{\prod_{j=2}^{4}(q^{j}-1)+1+q(K(\lambda ;
\beta^{-1})^{2}-q)\} (cf. \;(\ref{m}))
\end{equation}
\begin{equation*}
=q^{6}\{q^{2}(q^{2}-1)(q^{4}-q-1)+K(\lambda ; \beta ^{-1})^{2}\}.
\qquad \qquad \;
\end{equation*}
Now, invoking Theorem {\ref{N}}, we obtain the result.\\
\end{proof}
\end{corollary}

We are now ready to prove Theorem \ref{A}, which is the main result
of this paper. To do that, we apply Pless power moment identity in
(\ref{l}), with $B=C^{\bot}(SL(n,q))$. Then, in view of Proposition
\ref{G} and utilizing (\ref{f}), the left hand side of (\ref{l}) is
given by

\begin{align*}
\begin{split}
 &\sum_{a \in \mathbb{F}_{q}^{*}} w(c(a))^{h}\\
 &=\frac{1}{2^{h}}\sum_{a \in \mathbb{F}_{q}^{*}} (N-q^{{\binom {n}{2}}} K_{n-1}(\lambda ; a))^h\\
 &=\frac{1}{2^{h}}\sum_{i=0}^{h}(-1)^{i}{\binom {h}{i}} N^{h-i}q^{{\binom {n}{2}}i} MK_{n-1}^{i}\\
 &=\frac{1}{2^{h}} (-1)^{h}q^{{\binom{n}{2}}h} MK_{n-1}^{h}\\
   &+\frac{1}{2^{h}}\sum_{i=0}^{h-1} (-1)^{i}{\binom{h}{i}}N^{h-i}q^{{\binom {n}{2}}i} MK_{n-1}^{i}.\\
\end{split}
\end{align*}

On the other hand, noting that $ dim$ $C^{\bot}(SL(n,q))=r$ (cf.
Proposition \ref{G}), the right hand side of (\ref{l}) is given by

\begin{equation*}
q\sum_{i=0}^{min \{N,h\}}(-1)^{i}C_{i} \sum_{t=i}^{h}t!
S(h,t)2^{-t}{\binom {N-i}{N-t}}.
\end{equation*}

Here the frequencies $C_{i}$ of codewords with weight $i$ in $
C=C(SL(n,q))$ are given by (\ref{n}).

Now, Corollary \ref{B} follows from Theorem \ref{A} and Corollary
\ref{Q}, and Corollary \ref{S} from Theorem \ref{A} and Corollary \ref{R}.\\

\begin{corollary}\label{S}
For all positive integers $h$ , we have the following recursive
formula for the moments of the 3-dimensional Kloosterman sums
$MK_{3}^{h}$,
\begin{align*}
\begin{split}
q^{6h}MK_{3}^{h}&= \sum_{i=0}^{h-1} (-1)^{h+i+1}{\binom {h}{i}} N^{h-i}q^{6i}MK_{3}^{i}\\
 &+q \sum_{i=0}^{min \{N,h\}} (-1)^{h+i}C_{i}\sum_{t=i}^{h}t! S(h,t)2^{h-t}{\binom {N-i}{N-t}}.\\
\end{split}
\end{align*}

Here $ N=q^{6}\prod_{j=2}^{4}(q^{j}-1)$ is the order of $ SL(4,q)$,
$\{C_{i}\}_{i=0}^{N}$  denotes  the weight distribution of the code
$ C=C(SL(4,q))$ given by (\ref{o}), and $ S(h,t)$ indicates the
Stirling number of the second kind as in (\ref{b}).
\end{corollary}

%%%%%%%%%%%%%%%%%%%%%%%%%%%%%%%%%%%%%%%%%%%%%%%%%%%%%%%%%%%%%%%%%%%%
\section{Remarks}
%%%%%%%%%%%%%%%%%%%%%%%%%%%%%%%%%%%%%%%%%%%%%%%%%%%%%%%%%%%%%%%%%%%%

 Here we will briefly review the previous results on power moments of
Kloosterman sums $MK^{h}$, and make some comments on our result in
(\ref{c}). For any $q=p^{r}$ ($p$ a prime),\\

\begin{equation}\label{q}
MK^{h}=\frac{q^{2}}{q-1}A_{h}-(q-1)^{h-1}+2(-1)^{h-1},
\end{equation}
where
\begin{equation*} A_{h}=|\{(\alpha_{1},\cdots,\alpha_{h}) \in (\mathbb{F}_{q}^{*})^{h} | \;
\sum_{j=1}^{h}\alpha_{j}=0=\sum_{j=1}^{h}\alpha_{j}^{-1}\}|.
\end{equation*}

For $h \in \mathbb{Z}_{\geq 0}$ , define $M_h$  as:

\begin{align*}
M_{h}=|\{( \alpha_{1},\cdots,\alpha_{h})
\in(\mathbb{F}_{q}^{*})^{h}|\sum_{j=1}^{h}\alpha_{j}=1=\sum_{j=1}^{h}\alpha_{j}^{-1}
\}|,
\end{align*}
for $h>0$, and $M_{0}=0$.\\

Then, as one can see, $(q-1)M_{h-1}=A_h$, for any positive integer $h$. So (\ref{q}) can be rewritten as\\
\begin{equation}\label{r}
MK^h=q^2M_{h-1}-(q-1)^{h-1}+2(-1)^{h-1}(h \geq 1).
\end{equation}

Sali\'{e} obtained this form of expression for $MK^h$ already in
{\cite{HS}}, for any odd prime $q$. Iwaniec {\cite{HI}} showed the
expression (\ref{q}) for any prime $q$. However, the proof given
there works for any prime power $q$, without any restriction. Also,
this is a special case of Theorem \ref{A} in {\cite{HD}}, as
mentioned in Remark 2 there.

Let $q$ = $p$ be any prime. Then
\begin{align*}
 &MK^1=1,\\
 &MK^2=p^2-p-1,\\
 &MK^3=(\frac{-3}{p})p^2+2p+1\\
 &(\textmd{with the understanding }(\frac{-3}{2})=-1,(\frac{
 -3}{3})=0),
\end{align*}
\begin{equation*}
MK^4=
\begin{cases}
2p^3-3p^2-3p-1,& \text {if $p \geq 3$,}\\
1,& \text {if $p=2$}.
\end{cases}
\end{equation*}

Sali\'{e} obtained these results in {\cite{HS}} by determining
$M_1$, $M_2$, $M_3$, and Iwaniec got these ones in {\cite{HI}} by
computing $A_2$, $A_3$,$A_4$.

Except {\cite{L1}} for $1 \leq h \leq 4$ and {\cite{MK}} for $h=6$,
not much progress had been made until Moisio succeeded in evaluating
$MK^h$, for the other values of $h$ with $h \leq 10$ over the finite
fields of characteristic two (Similar results exist also over the
finite fields of characteristic three
{\cite{GR}},{\cite{M2}}). His results are as follows:\\

\begin{align*}
\begin{split}
MK^1=&1,\\
MK^2=&q^2-q-1,\\
MK^3=&(-1)^r q^2+2q+1,\\
MK^4=&2q^3-2q^2-3q-1,\\
MK^{5}=&(u_{1}+(-1)^{r}4)q^{3}+5q^{2}+4q+1,\\
MK^{6}=&5q^{4}-(5+(-1)^{r})q^{3}-9q^{2}-5q-1,\\
MK^{7}=&(u_{2}+6u_{1}+(-1)^{r}14+1)q^{4}+14q^{3}+14q^{2}+6q\\
       &+1,\\
MK^8=&14q^5-(15+(-1)^r7)q^4-28q^3-20q^2-7q-1,\\
MK^9=&(u_3+8u_2+27u_1+8+(-1)^r 48)q^5+42q^4+48q^3\\
     &+27q^2+8q+1,\\
\end{split}
\end{align*}
\begin{align}\label{s}
\begin{split}
MK^{10}=&42q^{6}-(51+(-1)^{r} 35)q^{5}-90q^{4}-75q^{3}-35q^{2}\\
        &-9q-1-u_{4}.\\
\end{split}
\end{align}

Here $u_1$, $u_2$, $u_3$, $u_4$ are the following numbers which are
dependent upon the extension degree $r$ of $\mathbb{F}_{q}$ over
$\mathbb{F}_{2}$:
\begin{align*}
\begin{split}
u_1=&((1+\sqrt{-15 })/4)^r+((1-\sqrt{-15})/4)^r,\\
u_2=&((-5+\sqrt{-39})/8)^r+((-5-\sqrt{ -39})/8)^r,\\
u_3=&((-3+\sqrt {505}+\sqrt{-510-6 \sqrt{505 }})/32)^r\\
    &+((-3+\sqrt{505}-\sqrt{{-510-6 \sqrt{505}}})/32)^r\\
    &+((-3-\sqrt{505}+\sqrt{-510+6\sqrt{505}})/32)^{r}\\
    &+((-3-\sqrt{505}-\sqrt{-510+6\sqrt{505}})/32)^{r},\\
u_4=&(-12+4\sqrt{-119})^r+(-12-4 \sqrt{-119})^r.\\
\end{split}
\end{align*}

As we mentioned earlier, these were obtained, via Pless power moment
identity, by expressing power moments of Kloosterman sums in terms
of the frequencies of weights in the binary Zetterberg code of
length $q+1$ . In fact, Moisio used the frequencies $B_{i}$ in the
Zetterberg code for $i \leq 12$, which were available in Table 6.2
of {\cite{RM}}.

Even though it was a breakthrough, it had a few drawbacks. Firstly,
the way it is proved is too indirect, since the frequencies are
expressed in terms of the Eichler Selberg trace formulas for the
Hecke operators acting on certain spaces of cusp forms for
$\Gamma_{1}(4)$. Secondly, the power moments of Kloosterman sums are
obtained only for $h \leq 10$ and not for any higher order moments.
On the other hand, our formula in (\ref{c}) allows one, at least in
principle, to compute moments of all orders for any given $q$.
Moreover, it gives a recursive formula not only for power moments of
Kloosterman sums but also for those of multi-dimensional Kloosterman
sums(cf. (\ref{a})). Nevertheless, obviously it is good to have
explicit formulas like the ones presented in ({\ref{s}}) . In the
next section, we will give some numerical examples demonstrating
that our formula in ({\ref{c}}) is quite useful for evaluating power
moments of Kloosterman sums for each given $q$ .

%%%%%%%%%%%%%%%%%%%%%%%%%%%%%%%%%%%%%%%%%%%%%%%%%%%%%%%%%%%%%%%%%%%%%%%
\section{Examples}
%%%%%%%%%%%%%%%%%%%%%%%%%%%%%%%%%%%%%%%%%%%%%%%%%%%%%%%%%%%%%%%%%%%%%%%
In this section, for small values of $i$, we compute, by using
Corollary \ref{B} and MAGMA, the frequencies $C_{i}$ of weights in
$C(SL(2,2^{3}))$ and $C(SL(2,2^{4}))$, and the power moments
$MK^{i}$ of Kloosterman sums over $\mathbb{F}_{2^{3}}$ and
$\mathbb{F}_{2^{4}}$. In particular, our results confirm those of
Moisio's given in (\ref{s}), when $q=2^{3}$ and $q=2^{4}$.

\begin{table}[!htp]
\begin{center}
\begin{tabular}{c c c c }
\multicolumn{4}{c}{TABLE I} \\
\multicolumn{4}{c}{The weight distribution of $C(SL(2,2^{3}))$} \\
\\
\hline
w & frequency & w& frequency\\[0.5pt]
\hline\\[0.5pt]
 \scriptsize{0} &\scriptsize{1}                    & \scriptsize{11}  & \scriptsize{1495424065262442956416}\\
 \scriptsize{1} &\scriptsize{64}                   & \scriptsize{12}  & \scriptsize{61437005346735099526740}\\
 \scriptsize{2} &\scriptsize{15844}                & \scriptsize{13}  &\scriptsize{ 2325154356197975713774208}\\
 \scriptsize{3} &\scriptsize{2650560}              & \scriptsize{14}  &\scriptsize{ 81546484920999191101202360}\\
 \scriptsize{4} &\scriptsize{332067914}            &\scriptsize{15}   &\scriptsize{ 2663851840752718923500482944}\\
 \scriptsize{5} &\scriptsize{33207770816}          & \scriptsize{16}  &\scriptsize{ 81413971883002952517354367429}\\
 \scriptsize{6} &\scriptsize{2761774095732}        & \scriptsize{17}  &\scriptsize{ 2337059898759141068388769445824}\\
 \scriptsize{7} &\scriptsize{196480443747136}      &\scriptsize{18}   & \scriptsize{63230453927539041393172170525052}\\
 \scriptsize{8} &\scriptsize{12206347634256355}    & \scriptsize{19}  &\scriptsize{ 1617368453093893435845237341156928}\\
 \scriptsize{9} &\scriptsize{672705382226871680}    & \scriptsize{20} & \scriptsize{39221184987526914436447793737809822}\\
 \scriptsize{10}&\scriptsize{33298916433035363704} &\scriptsize{21}   &\scriptsize{903954930188715550538753640492641088}\\
\hline
\end{tabular}
\end{center}
\end{table}

\begin{table}[!htp]
\begin{center}
\begin{tabular}{c c c c c c }
\multicolumn{6}{c}{TABLE II} \\
\multicolumn{6}{c}{The power moments of Kloosterman sums over
$\mathbb{F}_{2^{3}}$ }\\
\\
\hline
\\
$i$ & $MK^{i}$ & $i$ & $MK^{i}$ & $i$ & $MK^{i}$ \\[0.5pt]
\hline\\[0.5pt]
 \scriptsize{0}& \scriptsize{7}        &  \scriptsize{10} & \scriptsize{9942775}         & \scriptsize{20}  & \scriptsize{95377891993831}\\
 \scriptsize{1}& \scriptsize{1}        &  \scriptsize{11} & \scriptsize{-48296687}       & \scriptsize{21}  & \scriptsize{-476805777143519}\\
 \scriptsize{2}& \scriptsize{55}       &  \scriptsize{12} & \scriptsize{245734951}       & \scriptsize{22}  & \scriptsize{2384279934194455}\\
 \scriptsize{3}& \scriptsize{-47}      &  \scriptsize{13} & \scriptsize{-1215920159}     & \scriptsize{23}  & \scriptsize{-11920646525541647}\\
 \scriptsize{4}& \scriptsize{871}      &  \scriptsize{14} & \scriptsize{6117864535}      & \scriptsize{24}  & \scriptsize{59605492064000071}\\
 \scriptsize{5}& \scriptsize{-2399}    &  \scriptsize{15} & \scriptsize{-30474531407}    & \scriptsize{25}  & \scriptsize{-298020682011124799}\\
 \scriptsize{6}& \scriptsize{17815}    &  \scriptsize{16} & \scriptsize{152717030791}    & \scriptsize{26}  &  \scriptsize{1490123744982250615}\\
 \scriptsize{7}& \scriptsize{-71567}   &  \scriptsize{17} & \scriptsize{-762552032639}   & \scriptsize{27}  & \scriptsize{-7450557720131373167}\\
 \scriptsize{8}& \scriptsize{410311}   &  \scriptsize{18} & \scriptsize{3815859527095}   & \scriptsize{28}  &  \scriptsize{37252971614996505511}\\
 \scriptsize{9}& \scriptsize{-1894079}  & \scriptsize{19} & \scriptsize{-19069999543727} & \scriptsize{29}  & \scriptsize{-186264309031963608479}\\
\hline
\end{tabular}
\end{center}
\end{table}

\begin{table}[!htp]
\begin{center}
\begin{tabular}{c c c c }
\multicolumn{4}{c}{TABLE III} \\
\multicolumn{4}{c}{The weight distribution of $C(SL(2,2^{4}))$} \\
\\
\hline
w & frequency & w& frequency\\[0.5pt]
\hline\\[0.5pt]
 \scriptsize{0}  & \scriptsize{1}               & \scriptsize{6}  & \scriptsize{398943240589827320}\\
 \scriptsize{1}  & \scriptsize{256}             & \scriptsize{7}  & \scriptsize{232184965775802188544}\\
 \scriptsize{2}  & \scriptsize{520072}          & \scriptsize{8}  & \scriptsize{118211170698394115200330}\\
 \scriptsize{3}  & \scriptsize{706962176}       & \scriptsize{9}  & \scriptsize{53483987453818691622983424}\\
 \scriptsize{4}  & \scriptsize{720560061732}    & \scriptsize{10} & \scriptsize{21773331292449548118228026776}\\
 \scriptsize{5}  & \scriptsize{587401078798592} & \scriptsize{11} & \scriptsize{8056132578206330016084726166784}\\
\hline
\end{tabular}
\end{center}
\end{table}

\begin{table}[!htp]
\begin{center}
\begin{tabular}{c c c c c c }
\multicolumn{6}{c}{TABLE IV} \\
\multicolumn{6}{c}{The power moments of Kloosterman sums over
$\mathbb{F}_{2^{4}}$}\\
\\
\hline
\\
$i$ & $MK^{i}$ & $i$ & $MK^{i}$ & $i$ & $MK^{i}$ \\[0.5pt]
\hline\\[0.5pt]
 \scriptsize{0}  &\scriptsize{15}      & \scriptsize{4}   &  \scriptsize{7631}    &   \scriptsize{8}   & \scriptsize{13118351}\\
 \scriptsize{1}  & \scriptsize{1}      & \scriptsize{5}   &  \scriptsize{22081}   &   \scriptsize{9}   &\scriptsize{72973441}\\
 \scriptsize{2}  & \scriptsize{239}    & \scriptsize{6}   &   \scriptsize{300719} &   \scriptsize{10}  & \scriptsize{604249199}\\
 \scriptsize{3}  & \scriptsize{289}    & \scriptsize{7}   &   \scriptsize{1343329}&   \scriptsize{11}  & \scriptsize{3760049569}\\
\hline
\end{tabular}
\end{center}
\end{table}

\begin{center}
ACKNOWLEDGMENT
\end{center}

I would like to thank Mr. Dong Chan Kim for providing me with the
above tables.

%%%%%%%%%%%%%%%%%%%%%%%%%%%%%%%%%%%%%%%%%%%%%%%%%%%%%%%%%%%%%%%%%%%%%%%

%\begin{biography}[{\includegraphics[width=1in,height=1.25in,clip,keepaspectratio]{mshell}}]{Michael Shell}
% where an .eps filename suffix will be assumed under latex, and a .pdf suffix
% will be assumed for pdflatex; or if you just want to reserve a space for
% a photo:

%\begin{biography}{Michael Shell}
%\end{biography}

%\begin{biographynophoto}{John Doe}
%\end{biographynophoto}

% insert where needed to balance the two columns on the last page
%\newpage

%\begin{biographynophoto}{Jane Doe}
%\end{biographynophoto}


%\vfill



\begin{thebibliography}{20}
%%%%%%%%%%%%%%%%%%%%%%%%%%%%%%%%%%%%%%%%%%%%%%%%%%%%%%%%%%%%%%%%%%%%%%%
\bibitem{L1}
L. Carlitz,
\newblock {``Gauss sums over finite fields of order $2^{n}$,''}
\newblock { Acta Arith.,vol. 15 , pp. 247-265, 1969.}

\bibitem{L2}
 L. Carlitz,
\newblock {``A note on exponential sums,''}
\newblock { Pacific J. Math.,vol. 30 , pp. 35-37, 1969.}

\bibitem{HD}
Hi-joon Chae and D. S. Kim,
\newblock {``A generalization of power moments of Kloosterman sums,''}
\newblock { Arch. Math.(Basel), vol. 89 , pp. 152-156, 2007.}


\bibitem{GR}
G. van der Geer, R. Schoof and M. van der Vlugt,
\newblock {``Weight formulas for ternary Melas codes,''}
\newblock { Math. Comp., vol. 58 , pp. 781-792, 1992.}

\bibitem{HI}
H. Iwaniec,
\newblock {\emph{Topics in Classical Automorphic Forms, }}
\newblock {Amer. Math. Soc., Providence, R. I., 1997.}

\bibitem{DS}
D. S. Kim,
\newblock {``Gauss sums for general and special linear groups over a finite field,'' }
\newblock { Arch. Math.(Basel), vol. 69 , pp. 297-304, 1997.}

\bibitem{HDK}
H. D. Kloosterman,
\newblock {$\emph{``On the representation of numbers in the form }$
$ax^2+by^{2}+cz^{2}+dt^{2}$,''}
\newblock { Acta. Math. vol. 49 , pp. 407-464, 1926.}


\bibitem{GJ}
G. Lachaud and J. Wolfmann,
\newblock {``The weights of the orthogonals of the extended quadratic binary Goppa codes,''}
\newblock {  IEEE Trans. Inform. Theory,  vol. 36 , pp. 686-692, 1990.}

\bibitem{RH}
R. Lidl and H. Niederreiter,
\newblock {\emph{Finite Fields, 2nd ed.}}
\newblock {Cambridge, U. K.:Cambridge University Pless, 1997, vol. 20, Encyclopedia of Mathematics and Its Applications.}

\bibitem{FN}
 F. J. MacWilliams and N. J. A. Sloane,
\newblock {\emph{The Theory of Error Correcting Codes.}}
\newblock {Amsterdam, The Netherlands: North-Holland, 1998.}

\bibitem{M1}
M. Moisio,
\newblock {``The moments of a Kloosterman sum and the weight distribution of a Zetterberg-type binary cyclic code,''}
\newblock {IEEE Trans. Inform. Theory, vol. 53, pp. 843-847, 2007.}

\bibitem{M2}
M. Moisio,
\newblock {``On the moments of Kloosterman sums and fibre products of Kloosterman curves,''}
\newblock {Finite Field Appl., vol.14, pp. 515-531, 2008.}



\bibitem{M3}
M. Moisio,
\newblock {``Kloosterman sums, elliptic curves, and irreducible polynomials with prescribed trace and norm,''}
\newblock {\emph{Acta Arith., to appear.}}



\bibitem{MK}
M. Moisio and K. Ranto,
\newblock {``Klooserman sum identities and low-weight codewords in a cyclic code with two zeros,''}
\newblock {Finite Fields Appl.,vol.13, pp. 922-935, 2007.}

\bibitem{HS}
H. Sali\'{e},
\newblock {``Uber die Kloostermanschen Summen $S(u,v;q)$,''}
\newblock {Math. Z., vol. 34, pp. 91-109, 1931.}

\bibitem{RM}
R. Schoof and M. van der Vlugt,
\newblock {``Hecke operators and the weight distributions of certain codes,''}
\newblock {J. Combin. Theory Ser. A, vol. 57, pp.163-186, 1991.}

\end{thebibliography}
\end{document}